\documentclass[a4paper,12pt,reqno]{amsart}

\usepackage[letterpaper, hmargin=1in, top=1in, bottom=1.2in, footskip=0.7in]{geometry}
\usepackage{rotating,pdflscape,xcolor,amssymb,subfigure,psfrag,amsmath,eufrak,bbm,epsfig,amsthm,mathtools,fancyhdr,graphicx}
\definecolor{vdarkred}{rgb}{0.6,0,0.2}
\definecolor{vdarkblue}{rgb}{0,0.2,0.6}
\usepackage[pdftex, colorlinks, linkcolor=vdarkblue,citecolor=vdarkred]{hyperref}
\usepackage{a4wide,amscd}
\usepackage{changepage}

\usepackage{tikz} 
 \usepackage[usenames,dvipsnames]{pstricks}
 \usepackage{pst-grad} 
 \usepackage{pst-plot} 

\usepackage[small]{caption}

\newcommand{\ii}{\mathrm{i}}
\newcommand{\me}{\mathrm{e}}

\newcommand{\cC}{\mathcal{C}}

\newcommand{\bE}{\mathbf{E}}
\newcommand{\bG}{\mathbf{G}}

\newcommand{\bX}{\mathbf{X}}
\newcommand{\bY}{\mathbf{Y}}

\newcommand{\bee}{\mathbf{e}}

\newcommand{\bu}{\mathbf{u}}

\newcommand{\lam}{\lambda}

\newcommand{\al}{\alpha}

\newcommand{\Si}{\Sigma}

\newcommand{\del}{\delta}

\newcommand{\ld}{\ldots}

\newcommand{\beg}{\begin}
\newcommand{\en}{\end}

\renewcommand{\Im}{\mathfrak{Im}}
\renewcommand{\Re}{\mathfrak{Re}}

\newcommand{\trm}{\textrm}

\newcommand{\bgt}{\begin{itemize}}
\newcommand{\ent}{\end{itemize}}

\newcommand{\eqre}{\eqref}
\newcommand{\re}{\ref}
\newcommand{\la}{\label}

\newcommand{\brem}{\begin{rmk}}
\newcommand{\erem}{\end{rmk}}
\newcommand{\blem}{\begin{lem}}
\newcommand{\elem}{\end{lem}}
\newcommand{\bcor}{\begin{cor}}
\newcommand{\ecor}{\end{cor}}
\newcommand{\bTh}{\begin{Th}}
\newcommand{\eTh}{\end{Th}}
\newcommand{\bpropo}{\begin{propo}}
\newcommand{\epropo}{\end{propo}}

\newcommand{\op}{\operatorname}

\newcommand{\Tr}{\operatorname{Tr}}

\newcommand{\ud}{\mathrm{d}}

\newcommand{\E}{\op{\mathbb{E}}}

\newcommand{\R}{\mathbb{R}}
\newcommand{\C}{\mathbb{C}}

\newcommand{\p}{\mathbb{P}}

\newcommand{\pro}{probability }

\newcommand{\f}{\frac}
\newcommand{\ff}{\frac{1}}
\newcommand{\lf}{\left}
\newcommand{\ri}{\right}

\newcommand{\st}{such that }

\newcommand{\ti}{\times}

\newcommand{\mc}{\mathcal}

\newcommand{\eps}{\varepsilon}

\newcommand{\bck}{\backslash}

\newcommand{\bbm}{\begin{bmatrix}}
\newcommand{\ebm}{\end{bmatrix}}
\newcommand{\bes}{\begin{equation*}}
\newcommand{\ees}{\end{equation*}}
\newcommand{\be}{\begin{equation}}
\newcommand{\ee}{\end{equation}}
\newcommand{\beqy}{\begin{eqnarray}}
\newcommand{\eeqy}{\end{eqnarray}}
\newcommand{\beq}{\begin{eqnarray*}}
\newcommand{\eeq}{\end{eqnarray*}}

\newcommand{\ie}{i.e. }

\newcommand{\bpm}{\begin{pmatrix}}
\newcommand{\epm}{\end{pmatrix}}

\newcommand{\cd}{\cdots}

\newcommand{\hlam}{\widehat{\lam}}

\newcommand{\bpr}{\beg{proof}}
\newcommand{\epr}{\en{proof}}

\newcommand{\pa}{\partial}

\newcommand{\AND}{\qquad\trm{ and }\qquad}

\newcommand{\hbE}{\widehat{\bE}}

\newcommand{\co}[1]{Corollary \re{#1}}

\allowdisplaybreaks
\newtheorem{Th}{Theorem}[section]

\newtheorem{propo}[Th]{Proposition}

\newtheorem{lem}[Th]{Lemma}

\newtheorem{cor}[Th]{Corollary}

\theoremstyle{definition}
\newtheorem{rmk}[Th]{Remark}

\long\def\symbolfootnote[#1]#2{\begingroup
\def\thefootnote{\fnsymbol{footnote}}\footnote[#1]{#2}\endgroup}

\author{Florent Benaych-Georges} \address{CFM, 23 rue de l'Universit\'e, 75007 Paris, France}
\email{florent.benaych@gmail.com}
\date{\today}
\title{A short proof of  Ledoit-P\'ech\'e's RIE formula for  covariance  matrices}
\begin{document}
\maketitle
\beg{abstract}This is a short proof of Ledoit-P\'ech\'e's RIE formula for  covariance  matrices. The proof is based on the Stein formula, which gives a very simple way to derive the result. One of the advantages of this approach is that it shows that the only really needed hypothesis, for the machinery to work, is that the mean of the eigenvalues of the true covariance matrix and the largest of them have the same order.
\en{abstract}

 
 {\bf Notation:} For $M$ a real matrix, $M'$ denotes the transpose of $M$, $\|M\|_{\text{F}}$ denotes the Frobenius norm of $M$ defined at \eqre{opt_problem00} and $\|M\|$ denotes the operator norm of $M$ (with respect to the canonical Euclidian norms). For $Z$ a random variable,   $\E Z$   denotes the expectation of $Z$.

 \section{Model}
 Let    $n$ and let  $X\in \R^n $ be  a  centered Gaussian  (column) random vector with covariance $\Si$.
 Let $$\bX=\bpm X(1)&\cd &X(T)\epm\in \R^{n\ti T}$$ be a collection of independent copies of $X$. 
 The problem is to estimate the true covariance $\Si$ out of the \emph{empirical estimator}  $$\bE=\ff T\bX\bX'.$$
We are specifically looking for a Rotationally Invariant Estimator\footnote{Such estimators arise naturally in the Bayesian framework where $\Si$ has been chosen at random with a distribution of which we only know that it is invariant under the action of the orthogonal group by conjugation.}, \ie an estimator of $\Si$ which has the same eigenvectors as $\bE$ (thus differs from it only by the eigenvalues). 
We want this estimator to be \emph{optimal}, \ie to be solution of  \be\la{intro_optimal_eq}\op{argmin}_{\op{estimators}}
   \| \op{Estimator}-\Si\|_{\text{F}}\ee
   among the estimators   whose eigenvectors 
  are those of  the  empirical estimator $\bE$. Here,   $\|\cdot\|_{\text{F}}$  denotes the \emph{Frobenius norm} ,  \ie   the standard Euclidean norm on matrices:   \be\la{opt_problem00}\|M\|_{\text{F}}:=\sqrt{\Tr MM'}.\ee
 
 Let $\lam_1,\ld, \lam_n$ denote the eigenvalues of $\bE$, with associated eigenvectors $\bu_1, \ld, \bu_n$. 
 It is obvious, by invariance of the Frobenius norm under the action of the orthogonal group by conjugation, that the solution of \eqre{intro_optimal_eq} is the real symmetric matrix with the eigenvectors $\bu_1, \ld, \bu_n$   and respective associated eigenvalues $\bu_k'\Si\bu_k$ ($1\le k\le n$).
 
   It follows that the optimal RIE is the one given by the formula \be\la{hlmak}\hbE:=\sum_{k=1}^n\hlam_k\bu_k\bu_k'\quad \text{for }\hlam_k:=\bu_k'\Si\bu_k.\ee

   Introducing the measures $$\sum_{k=1}^n \del_{\lam_k}$$ and $$\sum_{k=1}^n \bu_k'\Si\bu_k\del_{\lam_k},$$   we can express the $\hlam_k$ of \eqre{hlmak} as Radon-Nikodym derivatives, and, using \eqre{muoutofg}, we have that 
     for any $\eps>0$ \st $[\lam_k-\eps, \lam_k+\eps]\cap\{\lam_1, \ld, \lam_n\}=\{\lam_k\}$, \be\la{0407184IO00AA}\hlam_k=\lim_{\eta\to 0}\f{\int_{\lam_k-\eps}^{\lam_k+\eps}\Im L(x+\ii\eta)\ud x}{\int_{\lam_k-\eps}^{\lam_k+\eps}\Im G(x+\ii\eta)\ud x},\ee where $L$ and $G$ are the meromorphic functions on the complex plane   defined by $$L(z):=\ff T\Tr \bG\Si,\qquad G(z):=\ff T\Tr \bG,$$ for $\bG=\bG(z):=(z-\bE)^{-1}$
   the resolvent matrix  of $\bE$. 
   
   The main results of \cite{RIE0} can be summed-up in the following one (with main consequence Equation \eqre{mainformula2023} below):
   \bTh\la{mainth} For any fixed $z\in \C\bck\R$,  we have \be\la{mainformula}L(z)=1-\ff{1-q+zG(z)}+o(1),\ee where $o(1)$ denotes a random variable tending to zero in \pro as $n,T$ tend to infinity in such a way that $q:=n/T$ stays bounded as well as the operator norm of $\Si$ and $T/\Tr \Si$.\eTh

   \brem Note that  as $z$ tends to $\lam$ in \eqre{mainformula},  
 $$
 \f{\Im L(z)}{\Im G(z)}\sim \f{-\Im  \ff{1-q+zG(z)}}{\Im G(z)} 
 =  
 \ff{|1-q+z G(z)|^2}\f{\Im(zG(z))}{\Im G(z)}
 \sim \f\lam{|1-q+\lam G(z)|^2},
 $$ 
 which, given $G$ is normalized by $T$ and not by $n$, is the well known  RIE formula one can find e.g. in \cite{BBPo}. Note that in the finite dimensional case where this formula is meant to be applied, $G(z)$ is singular at any eigenvalue $\lam$, hence one has to slightly regularize the denominator, \ie use $G(z)$ for $z=\lam+\ii\eta$ with a positive spectral resolution $\eta=T^{-\al}$ (any exponent $\al \in (0, 1)$, typically $\al=1/2$,  should work). To sum up, for such a spectral resolution $\eta$, the formula of the cleaned eigenvalue $\hlam$ associated to an eigenvalue $\lam$ is given by \be\la{mainformula2023}\hlam=\f\lam{|1-q+\lam G(\lam+\ii\eta)|^2}.\ee
    \erem
   
   \brem One of the advantages of this approach is that it shows that the only really needed hypothesis, for the machinery to work, is that the mean $\Tr \Si/n$ of the eigenvalues of the true covariance matrix  and the largest of them (i.e. the operator norm $\|\Si\|$ of $\Si$) have the same order.
   \erem 
   
Let us now show how concentration of measure and the Stein formula for Gaussian random vectors allow to recover the above result very easily. In the case where the random vector $\bX$ is not Gaussian, what follows also works with a bit more of work, using e.g.  \cite[Lem. 1.13.9]{FloAntti} instead of the Stein formula and the McDiarmid inequality instead of Gaussian measure concentration. 
The idea of the proof is the following one: one starts with the formula $z\bG-I_n=\bG\bE$, pass to the expectation and expand  the RHT  using Stein formula for Gaussian vectors. The expansion makes products of traces appear, each of which concentrates around its expectation by the concentration of measure principle. At the end, we obtain the desired  relation between $G(z)$ and $L(z)$.

\bpr
 Let $$H:=\ff T\Tr \bG\bE.$$
 Using \co{cor91171}, 
 \beq \E \Tr \bG\bE&=&\E \ff T\sum_t \Tr \bG X(t)X(t)'\\
 &=&  \ff T\sum_t \E X(t)' \bG X(t)\\
 &=&\ff T\sum_t \E \Tr \bG \Si+\ff T\sum_t\sum_{k=1}^n \E \bee_k'\Si\lf(\f{\pa}{\pa X(t)_k}\bG \ri)X(t)
 \eeq
Note that  for any $t$, at $X(1), \ld, \widetilde{X(t)}, \ld X(T)$ fixed, 
the differential, at $X(t)\in\R^n$, of the function $X(t)\mapsto \bG$ is the function 
$$x\in \R^n\mapsto \ff T\bG\lf(X(t)x'+xX(t)'\ri)\bG,
$$
so that
$$\f{\pa}{\pa X(t)_k}\bG=\ff T\bG\lf(X(t)\bee_k'+\bee_kX(t)'\ri)\bG$$ and 
 \beq \E \Tr \bG\bE&=&\ff T\sum_t \E \Tr \bG \Si+\ff{T^2}\sum_t\sum_{k=1}^n \E \bee_k'\Si\lf(\bG\lf(X(t)\bee_k'+\bee_kX(t)'\ri)\bG\ri)X(t)\\
 &=& \E \Tr \bG \Si+\ff{T^2}\sum_t\sum_{k=1}^n \E \lf(\bee_k'\Si\bG X(t)\bee_k'\bG X(t)+\bee_k'\Si\bG  \bee_kX(t)'\bG X(t)\ri)\\
  &=& \E \Tr \bG \Si+\ff{T^2}\sum_t\E \lf(X(t)'\bG \Si\bG X(t) +\Tr( \Si\bG  )X(t)'\bG X(t)\ri)\\
  &=&  \E \Tr \bG \Si+\ff{T}\E \Tr \bG \Si\bG \bE +\ff{T}\E\Tr( \Si\bG  )\Tr (\bG \bE)
 \eeq
 Dividing by $T$ and using Proposition \re{247153} with \cite[Lem. B.2]{FloRomain1}, we find 
\be\la{28112033}H=L+LH+o(1).\ee 
Now, we are close to the conclusion but have to divide by $1+H$. 
 Writing $\bG\ti (z-\bE)=I_n$, we have $$z\bG-I_n=\bG\bE,$$ so that   $$H=zG-q=\f{q}{n}\sum_{i=1}^n\f{\lam_i}{z-\lam_i}.$$  
 Note that for $x=\Re z$, $\eta=\Im z$, we have \be\la{eqImH}|\Im H| =\f{q|\eta|}{n}\sum_{i=1}^n\f{\lam_i}{(x-\lam_i)^2+\eta^2}\ge \f{|\eta| \Tr\bE/T}{2x^2+2\|\bE\|^2+\eta^2}.\ee
By \eqre{eqImH} and Lemma \re{bEbounded}, there is $c>0$ \st with \pro tending to $1$, $|1+H|>c$. 

Hence by \eqre{28112033},  $$L=\f{H}{1+H}+o(1)=1-\ff{1+H}+o(1).$$\epr

 \blem\la{bEbounded}  As $n,T$ tend to infinity in such a way that $q=n/T$ stays bounded and $\Si$ stays bounded in operator norm, there is a constant $C>0$ (depending on the bounds on $q$ and on $\|\Si\|$) \st with \pro tending to $1$, the operator norm $\|\bE\|$ of $\bE$  and its trace $\Tr \bE$  satisfy $$\|\bE\|\le q\|\Si\|C\AND \Tr \bE\ge \Tr \Si/C$$ \elem
 
 \bpr Let $\bY\in \R^{n\ti T}$ be a matrix with independent standard Gaussian entries. Then $\bX$ has the same distribution as $\Si^{1/2}\bY$ and $T\|\bE\|$ has the same distribution as $\|\Si\bY\bY'\|$. Thus to prove the part about the operator norm,  it suffices to prove that there is a universal constant $C_0$ \st   with \pro tending to $1$, $\|\bY\bY'\|\le (n+T)C_0$. This is a very well known fact, that can be obtained for example as follows: given the eigenvalues of $\bY\bY'$ are, up to the null eigenvalue, the squares of the eigenvalues of $\bpm 0&\bY\\ \bY'&0\epm$, it follows from \cite[Th. 1.13.17]{FloAntti}. The part about the trace  follows from the fact that $\Tr \bE$ has the same law as $T^{-1}\Tr \Si \bY\bY'$ which has expectation $\Tr \Si$ and variance $\le 2T^{-1} \Tr \Si^2$.
  \epr
       \section{Appendix}\subsection{Stieltjes transform inversion}Any 
    signed measure $\mu$ on $\R$ can be recovered out of its Stieltjes transform\be\la{def:Stieltjes}g_\mu(z):=\int\f{\ud\mu(t)}{z-t} , \quad z\in \C\bck\R\ee by the formula   \be\la{muoutofg}\mu=-\ff\pi\lim_{\eta\to 0^+}(\Im g_\mu(x+\ii\eta)\ud x),\ee where the limit holds in the weak topology   (see e.g. \cite[Th. 2.4.3]{agz} and use the decomposition of any signed measure as a difference of finite  positive measures). 
    
%

 \subsection{Stein formula for   real Gaussian random vectors}
 \beg{propo}\la{SteinMultidim}Let $X=(X_1, \ld, X_d)$ be a real centered Gaussian vector with covariance $\Si$ and $f : \R^d  \to \R$ be a $\cC^1$  function with gradient having at most polynomial growth at infinity. Then for all $i_0=1, \ld, d$, $$\E X_{i_0}f(X_1, \ld, X_d) \;=\;\sum_{k=1}^d\Si_{i_0k}  \E (\partial_kf)(X_1, \ld, X_d). $$  
 \en{propo}

 \bcor\la{cor91171} With the same notation, considering $X$ as a column vector,  for $F:\R^d\to\R^{d\ti d}$ a matrix-valued function with gradient having at most polynomial growth at infinity. 
 we have   \beqy\la{2307181} \E X'F(X)X &=&\Tr \Si\E F(X)+\sum_{k=1}^d\lf(\E \Si (\pa_kF)(X) X\ri)_k .  
\eeqy  \ecor
 
 \bpr 
 We have, by  Proposition \re{SteinMultidim}, \beq \E X'F(X)X &=&\sum_{ij}\E X_{i}X_{j}F(X)_{ij}\\
  &=&\sum_{ijk}\E \Si_{ik}\f{\pa}{\pa X_k}X_{j}F(X)_{ij}\\
    &=&\sum_{ijk}\E \Si_{ik} \lf(\del_{j=k}F(X)_{ij}+X_{j}(\pa_kF)(X)_{ij}\ri)\\
     &=&\Tr \Si\E F(X)+\sum_{ijk}\E \Si_{ik}  X_{j}(\pa_kF)(X)_{ij} \\
          &=&\Tr \Si\E F(X)+\sum_{k}\lf(\E \Si (\pa_kF)(X) X\ri)_k
 \eeq
 \epr
 
  \subsection{Concentration of measure for Gaussian vectors}The following proposition can be found   e.g. in   \cite[Sec. 4.4.1]{agz} or \cite[Th. 5.2.2]{Vershynin}.
 \beg{propo}\la{247153}Let $X=(X_1, \ld, X_d)$ be a standard real Gaussian vector and $f : \R^d  \to \R$ be a $\mc{C}^1$  function with  gradient $\nabla f$. Then we have \be\la{247151}\op{Var}(f(X))\;\le\; \E \|\nabla f(X)\|^2,\ee where $\|\,\cdot\,\|$ denotes the standard Euclidean norm. 

Besides, if $f$ is $k$-Lipschitz, then for any $t>0$, we have   \be\la{247152}\p (|f(X)-\E f(X)|\ge t)\; \le \; 2\me^{-\f{t^2}{2k^2}},
\ee \ie  $f(X)-\E f(X)$ is Sub-Gaussian with Sub-Gaussian norm $\le k$, up to a universal constant factor. 
 \en{propo}

  \begin{thebibliography}{10}
  \bibitem{agz} Anderson, G., Guionnet, A., Zeitouni, O.  \emph{An Introduction to Random Matrices}. Cambridge Studies in Advanced Mathematics, {118} (2009).
    \bibitem{FloRomain1} Benaych-Georges, F., Couillet, R. \emph{Spectral analysis of the Gram matrix of mixture models}, ESAIM Probab. Statist., Vol. 20 (2016), 217--237. 
     \bibitem{FloAntti}  Benaych-Georges, F., Knowles, A. \emph{Local semicircle law for Wigner matrices. Advanced topics in random matrices}, 1--90, Panor. Synth\`eses, 53, Soc. Math. France, Paris, 2017.
 \bibitem{BBPo} Bun, J., Bouchaud,  J.-P., Potters, M. \emph{Cleaning large correlation matrices: Tools from Random Matrix Theory}, Physics Reports
Volume 666, Review article,   1--109, 2017.
\bibitem{RIE0} Ledoit, O.,   P\'ech\'e, S. \emph{Eigenvectors of some large sample covariance matrix ensembles} Probability Theory and Related Fields, 2011, 151.1, 233--264
 \bibitem{Vershynin} Vershynin, R. \emph{High-Dimensional Probability: An Introduction with Applications in Data Science}, Cambridge Series in Statistical and Probabilistic Mathematics, 2018.
 \en{thebibliography}
\en{document}